\documentclass[12pt]{amsart}

\newtheorem{thm}{Theorem}

\newtheorem{lem}[thm]{Lemma}
\newtheorem{pro}[thm]{Proposition}
\newtheorem{thmofothers}{Theorem}

\newtheorem{rem}{Remark}

\begin{document}
\title[Extension to a Martingale Inequality]{An Extension to the
Strong Domination Martingale Inequality}
\author{Stephen Montgomery-Smith}
\address{Department of Mathematics\\
University of Missouri\\
Columbia, Missouri 65211, USA.}
\email{stephen@math.missouri.edu}
\urladdr{http://www.math.missouri.edu/\~{}stephen}
\author{ Shih-Chi Shen }
\address{Department of Mathematics\\
University of Missouri\\
Columbia, Missouri 65211, USA.}
\email{mathgr75@math.missouri.edu}
\thanks{Both authors were supported in part by the NSF
and the Research Board of the University of Missouri}

\begin{abstract}
For each $1<p<\infty$, there exists a positive constant $c_p$,
depending only on $p$, such
that the following holds.  Let $(d_k)$, $(e_k)$ be real-valued martingale
difference sequences. If for for all bounded nonnegative predictable
sequences $(s_k)$ and all positive integers $k$ we have
\[E[s_k \vee |e_k|]\leq E[s_k \vee |d_k|]\]
then for all positive integers $n$ we have
\[ \left\| \sum_{k=1}^n e_k \right\|_p
\leq  c_p \left\| \sum_{k=1}^n d_k \right\|_p .\]
\end{abstract}
\maketitle

\section{Introduction}
Let $(\Omega ,\mathcal{F},P)$ be a probability space, and let
$(\mathcal{F}_k)$ be a filtration on $(\Omega ,\mathcal{F},P)$.
(We will suppose that $\mathcal{F}_0 = \{\emptyset,\Omega\}$.) If
an adapted sequence $(d_k)$ is a real-valued martingale difference
sequence, Burkholder's inequality \cite{B2} shows that for any
$1<p<\infty$, if $(v_k)$ is a predictable sequence bounded in
absolute value by 1, then there exists a positive constant $c_p$,
depending only on $p$, such that such that for all positive
integers $n$
\[ \left \| \sum_{k=1}^n v_k d_k \right\|_p
\leq c_p \left\| \sum_{k=1}^n d_k \right\|_p .\]
Later Burkholder \cite{B3} extended this result to
\emph{subordination} martingales: if $(d_k)$, $(e_k)$ are two
martingale difference sequences such that $(e_k)$ is subordinate
to $(d_k)$, that is, for all $k \geq 1$,
\begin{equation}
\label{subordination}|e_k|\leq |d_k|
\end{equation}
then there exists a positive constant $c_p$, depending only on
$p$, such that for all positive integers $n$
\begin{equation}
\label{ek<=dk} \left\|\sum_{k=1}^n e_k \right\|_p \leq c_p \left\|
\sum_{k=1}^n d_k \right\|_p.
\end{equation}
A different approach to this inequality was proposed by Kwapie\'n
and Woyczi\'{n}ski \cite{KW1} (see also \cite{KW2}). Two adapted
sequence $(d_k)$ and $(e_k)$ are said to be \emph{tangent} if for
each $k\geq 1$, we have that the law of $d_k$ conditionally on $
\mathcal{F}_{k-1} $ is the same as the law of $e_k$ conditionally
on $\mathcal{F}_{k-1}$, that is,
\begin{equation}
\label{tangent}P(d_k > \lambda | \mathcal{F}_{k-1})= P(e_k >
\lambda | \mathcal{F}_{k-1})
\end{equation}
for all real numbers $\lambda $. Answering a conjecture of
Kwapie\'n and Woyczi\'nski \cite{KW1}, it was proved by Hitczenko
\cite{H1} and Zinn \cite{Z} that for $1<p<\infty$ that there
exists a positive constant $c_p$, depending only on $p$, such that
if $(d_k)$ and $(e_k)$ are martingale difference sequences and
$(d_k)$, $(e_k)$ are tangent, then for all positive integers $n$
we have equation~(\ref{ek<=dk}).

Given two adapted sequences, $(e_k)$ is said to be \emph{strongly
dominated by} $(d_k)$ if for each $k\geq 1$,
\begin{equation}
\label{strong-d}P(|e_k| > \lambda | \mathcal{F}_{k-1})\leq P(|d_k|
> \lambda | \mathcal{F}_{k-1})
\end{equation}
for all $\lambda \geq 0$. It is obvious that the case of
(\ref{subordination}) and the case of (\ref{tangent}) are
contained in the cases of (\ref{strong-d}). Thus the following result of
Kwapie\'n and
Woyczi\'nski \cite{KW1} is a common generalization of these
two results: if $(d_k)$, $(e_k)$ are two martingale difference
sequences such that $(e_k)$ is strongly dominated by $(d_k)$, then
there exists a positive constant $c_p$, depending only on $p$,
such that for all positive integers $n$ equation (\ref{ek<=dk})
holds.

The purpose of this paper is to use a different approach to
provide another common generalization of those two results, an
even a further extension to Kwapie\'n and Woyczi\'nski's
result.

\begin{thm}
\label{t ek<=dk sk}
For each $1<p<\infty$, there exists a positive constant $c_p$,
depending only on $p$, such
that the following holds.
Let $(d_k)$, $(e_k)$ be real-valued martingale
difference sequences. If for for all bounded nonnegative predictable
sequence $(s_k)$ and all positive integers $k$ we have
\begin{equation}
\label{skek<=skdk}
E[s_k\vee |e_k|]\leq E[s_k \vee |d_k|]
\end{equation}
then for all positive integers $n$ we have
equation~(\ref{ek<=dk}).
\end{thm}

\begin{rem} {\upshape
We have that (\ref{skek<=skdk}) is
equivalent to
\begin{equation}
\label{skekFk-1<=skdkFk-1} E[(\lambda \vee
|e_k|)|\mathcal{F}_{k-1}]\leq E[(\lambda \vee
|d_k|)|\mathcal{F}_{k-1}]
\end{equation}
for all $\lambda \geq 0$. This is because for any $A_k \in
\mathcal{F}_{k-1}$ and $a \geq 0$ we have that 
$(a \chi_{A_k^c} \vee \lambda)$ is
predictable, and hence
\[E[(a \chi_{A_k^c} \vee \lambda ) \vee |e_k| - a \chi_{A_k^c}]
\leq E[(a \chi_{A_k^c} \vee \lambda ) \vee |d_k| - a \chi_{A_k^c}]\]
When $a$ intends to infinity, we obtain
\[E[( \lambda \vee |e_k|)\chi_{A_k}]\leq E[(\lambda \vee |d_k|)\chi_{A_k}]\]
which is equivalent to (\ref{skekFk-1<=skdkFk-1}).
}\end{rem}

\begin{rem} {\upshape
To see  that Theorem \ref{t ek<=dk sk} is really an extension
to Kwapie\'n and Woyczi\'nski's result, we just simply observe
that (\ref{strong-d}) is equivalent to
$$P(\{|e_k|>\lambda \} \cap A_k) \leq P(\{|d_k|>\lambda \} \cap
A_k),$$ and (\ref{skekFk-1<=skdkFk-1}) is equivalent to
$$\int_\lambda^\infty P(\{|e_k|>t \} \cap A_k)dt \leq \int_\lambda^\infty P(\{|d_k|>t \} \cap A_k)dt$$
for all $A_k \in \mathcal{F}_{k-1}$.
}\end{rem}

\begin{rem} {\upshape
Once we have Theorem~\ref{t ek<=dk sk}, we can obtain that for
$\kappa \geq 1$, if $$P(|e_k| > \lambda | \mathcal{F}_{k-1})\leq
\kappa P(|d_k| > \lambda | \mathcal{F}_{k-1}),$$ we have
\begin{equation}
\label{kappa} \left\|\sum_{k=1}^n e_k \right\|_p \leq  \kappa c_p
\left\| \sum_{k=1}^n d_k \right\|_p.
\end{equation}
This is because 
\begin{eqnarray*}
\int_\lambda^\infty P(\{|e_k|>t \} \cap A_k)dt
&\leq& \kappa \int_\lambda^\infty P(\{|d_k|>t \} \cap A_k)dt\\
&=& \kappa \int_{\frac{\lambda}{\kappa}}^\infty P\left( \left\{|d_k|> \frac{t}{\kappa } \right\} \cap A_k \right)d \left(\frac{t}{\kappa}\right)\\
&\leq& \int_\lambda^\infty P(\{\kappa |d_k|>t \} \cap A_k)dt\\
\end{eqnarray*}
Hence $$E[(\lambda \vee |e_k|)|\mathcal{F}_{k-1}]\leq E[(\lambda
\vee \kappa |d_k|)|\mathcal{F}_{k-1}]$$ and equation (\ref{kappa})
follows.
}\end{rem}

Let us give an application of Theorem~\ref{t ek<=dk sk}.
In fact this application is essentially equivalent to Theorem~\ref{t ek<=dk sk},
and indeed will play a large role in its proof.
We will consider the probability space $[0,1]^\mathbb{N}$ equipped
with the product Lebesgue measure $\mathcal L$, and consider the
filtration $(\mathcal{L}_k)$, where $\mathcal L_k$ is the minimal
$\sigma$-field for which the first $k$ coordinate functions of
$[0,1]^{\mathbb{N}}$ are measurable. Then two sequences $(d_k)$
and $(e_k)$ are tangent if
\[ e_k (x_1,\dots,x_k)=d_k(x_1,\dots,x_{k-1},\phi_k(x_1,\dots,x_k))\]
where $(\phi_k:[0,1]^k\rightarrow [0,1])$ is a sequence of
measurable functions such that $\phi_k(x_1,\dots,x_{k-1},\cdot)$
is a measure preserving map for almost all $x_1$,\dots, $x_{k-1}$.

We will consider a more general situation.  Suppose we have a
sequence of linear operators $(T_k(x_1,\dots,x_{k-1}))$, depending
measurably upon $(x_k) \in [0,1]^{\mathbb{N}}$, that are bounded
operators on both $L_1([0,1])$ and $L_\infty ([0,1])$ with norm 1.
Then consider the condition
\begin{equation}
\label{ek=Tk(dk)} e_k(x_1,\dots,x_{k-1}, \cdot
)=[T_k(x_1,\dots,x_{k-1})]d_k(x_1,\dots,x_{k-1}, \cdot ).
\end{equation}

\begin{thm}
\label{t ek<=dk T}
For each $1<p<\infty$, there exists a positive constant $c_p$,
depending only on $p$, such
that the following holds.
If $(d_k)$, $(e_k)$ and $(T_k)$ are as above satisfying (\ref{ek=Tk(dk)}),
then for all positive integers $n$ we have
equation~(\ref{ek<=dk}).
\end{thm}

We will also need the following intermediate result.
For any random variable $f$,
let $f^{\#}$ be the decreasing rearrangement of $|f|$, that is,
\[f^{\#}(t)=\sup\{ s\in \mathbb{R}: P(|f|<s)<t\}.\]

\begin{thm}
\label{t ek<=dk sharp} For each $1<p<\infty$, there exists a
positive constant $c_p$, depending only on $p$, such that the
following holds. Let $(d_k)$, $(e_k)$ be martingale difference
sequences on $[0,1]^\mathbb{N}$ with respect to $(\mathcal{L}_k)$.
Suppose that for each positive integer $k$
\[\int_0^t (e_k(x_1,\dots,x_{k-1},\cdot))^{\#}(s)ds
\leq \int_0^t  (d_k(x_1,\dots,x_{k-1},\cdot))^{\#}(s)ds\] for all
$t\in [0,1]$ and almost all $x_1$,\dots,$x_{k-1}$. Then for all
positive integers $n$ we have equation~(\ref{ek<=dk}).
\end{thm}

\section{The Discrete Type Case}

In this section we will prove Theorems~\ref{t ek<=dk T} and~\ref{t
ek<=dk sharp} in a special discrete situation, which we now
describe. For any positive integer $N$, let $\Sigma_N $ be the
$\sigma $-field generated by the partition
$\{[\frac{i-1}{N},\frac{i}{N}): i=1,2,\dots,N\}$. Define a
filtration $(\mathcal{F}_k)$ on $[0,1]^{\mathbb{N}}$ by
$\mathcal{F}_k=\mathcal{L}_{k-1}\otimes \Sigma_N $. Suppose
$(d_k)$, $(e_k)$ are $(\mathcal{F}_k)$-adapted.  Then for each $k$
and for each $x_1$,\dots, $x_{k-1}$, we see that
$d_k(x_1,\dots,x_{k-1},\cdot)$ and $e_k(x_1,\dots,x_{k-1},\cdot)$
are $\Sigma_N$-measurable simple functions on $[0,1)$. Therefore
$d_k$ and $e_k$ can be written as $N$-dimensional vectors and
$T_k(x_1,\dots,x_{k-1})$ can be represented by a $N\times N$
matrix, that is,
\[
\left[ \begin{array}{c} e_k(1) \\ e_k(2) \\ \vdots \\ e_k(N)
\end{array} \right ]
= \left[\begin{array}{ccc} a_k(1,1), &\dots& ,a_k(1,N) \\
a_k(2,1), &\dots& ,a_k(2,N) \\ \vdots & & \vdots \\ a_k(N,1),
&\dots& ,a_k(N,N)\end{array}\right] \left[ \begin{array}{c} d_k(1)
\\ d_k(2) \\ \vdots \\ d_k(N) \end{array} \right]
\]
where
\[d_k(i)=d_k(x_1,\dots,x_{k-1},i)=d_k(x_1,\dots,x_{k})
\mbox{ if $x_k \in [\frac{i-1}{N},\frac{i}{N})$}\]
\[e_k(i)=e_k(x_1,\dots,x_{k-1},i)=e_k(x_1,\dots,x_{k})
\mbox{ if $x_k \in [\frac{i-1}{N},\frac{i}{N})$}\]
\[T_k=T_k(x_1,\dots,x_{k-1})=\left[(a_k(x_1,\dots,x_{k-1}))(i,j)\right]_{N\times N}=\left[a_k(i,j)\right]_{N\times N}\]
The condition of being martingale difference sequences implies that
\[ \sum_{i=1}^N d_k(i)= \sum_{i=1}^N e_k(i)=0\]

\begin{pro}
\label{t ek<=dk T disc}
Theorem~\ref{t ek<=dk T} holds in the case that
$(d_k)$ and $(e_k)$ are
adapted to the filtration $(\mathcal F_k)$ described above.
\end{pro}

In this discrete case, the boundedness of $\|T_k\|_{L_1([0,1])}$
and $\|T_k\|_{L_\infty ([0,1]) }$ by 1 is equivalent to the
condition that $\sum_{j=1}^N |a_k(i,j)|\leq 1 $ for all $i$ and
$\sum_{i=1}^N |a_k(i,j)|\leq 1$ for all $j$. We claim that without
loss of generality, we can assume that every row sum and column
sum of $T_k$ is 0, that is,
\[ \sum_{j=1}^N a_k(i,j)= \sum_{i=1}^N a_k(i,j)=0\]
for all $i$ and $j$.
Suppose the $i^{th}$ row sum $\sum_{j=1}^N
a_k(i,j)=R_k(i)$. Let $T'_k$ be the liner operator defined by
\[T'_k=\left[a_k(i,j)-\frac{R_k(i)}{N} \right]_{N\times N}\]
It is clear that every row sum of $T'_k$ is 0 and
\begin{eqnarray*}
(T'_kd_k)(i)&=& \sum_{j=1}^N \left(a_k(i,j)-\frac{R_k(i)}{N} \right)d_k(j)\\
            &=& \sum_{j=1}^N a_k(i,j)d_k(j)-\frac{R_k(i)}{N} \sum_{j=1}^N d_k(j)\\
            &=& e_k(i)
\end{eqnarray*}

Now we can assume that every row sum of $T_k$ is 0. Similarly
suppose the $j^{th}$ column sum $\sum_{i=1}^N a_k(i,j)=C_k(j)$.
Let $T''_k$ be the linear operator defined by
\[T''_k=\left[a_k(i,j)-\frac{C_k(j)}{N} \right]_{N\times N}\]
Again it is clear that every row sum and column sum of $T''_k$ is
0 and
\begin{eqnarray*}
(T''_kd_k)(i)&=& \sum_{j=1}^N \left(a_k(i,j)-\frac{C_k(j)}{N} \right)d_k(j)\\
            &=& \sum_{j=1}^N a_k(i,j)d_k(j)-\frac{1}{N} \sum_{j=1}^N C_k(j)d_k(j)\\
            &=& e_k(i)
\end{eqnarray*}
since
\[\sum_{i=1}^N e_k(i)= \sum_{j=1}^N C_k(j)d_k(j)=0\]
After adjusting $T_k$, it is easy to check that the norms of $T_k$
may be enlarged up to 4. Of course, we can pick up $T_k /4 $
instead and absorb the 4 into the constant $c_p$.

A nonnegative real matrix is said to be \emph{doubly stochastic}
if each of its row and column sum is 1. A sub-doubly stochastic
matrix means that each of its row and column sum is less than or
equal to 1. Therefore we can change the assumption in
Proposition~\ref{t ek<=dk T disc} to be that: ``for almost all
$x_1$,\dots, $x_{k-1}$, every row sum and column sum of the matrix
from $T_k$ is 0, and the matrix from $|T_k|$ is sub-doubly
stochastic for each positive integer $k$''

One of the fundamental results in the theory of doubly stochastic
matrices was introduced by Birkhoff \cite{BI1} (or see for example
\cite{M1}).

\begin{thmofothers}
\label{t double stoch}
If $M$ is a doubly stochastic matrix, then
\[ M=\sum_{i=1}^S \theta_iP_i\]
where $P_i$ are permutation matrices, and the $\theta _i$ are
nonnegative numbers satisfying $\sum_{i=1}^S \theta _i =1$.
\end{thmofothers}
\bigskip

\begin{lem}
\label{l 2nx2n}
If $M$ is a $n\times n$ sub-doubly stochastic matrix,
then there exists a
$2n\times 2n$ doubly stochastic matrix such that its upper left
$n\times n$ sub-matrix is $M$.
\end{lem}

\begin{proof} Suppose that $R(i)$ is the $i^{th}$ row sum of $M$,
$C(j)$ is the $j^{th}$ column sum and $S$ is the sum of all
entries. Let
\[A=\left[\begin{array}{ccc} \frac{1-R(1)}{n}, &\dots& ,\frac{1-R(1)}{n} \\
\vdots & & \vdots \\ \frac{1-R(n)}{n}, &\dots&
,\frac{1-R(n)}{n}\end{array}\right]_{n\times n}\]
\[B=\left[\begin{array}{ccc} \frac{1-C(1)}{n}, &\dots& ,\frac{1-C(n)}{n} \\
\vdots & & \vdots \\ \frac{1-C(1)}{n}, &\dots&
,\frac{1-C(n)}{n}\end{array}\right]_{n\times n}\]
\[C= \mbox{Diag}\left[\begin{array}{ccc} \frac{S}{n}, &\dots& ,\frac{S}{n}
\end{array}\right]_{n\times n}\]
Then define
\[M'=\left[\begin{array}{cc} M & A \\
B & C \end{array}\right]_{2n\times 2n}\] It is easy to check that
$M'$ is a doubly stochastic matrix.
\end{proof}

\begin{lem}
\label{l M+N}
If $M$ is a sub-doubly stochastic matrix, then there exists a sub-doubly
stochastic matrix $N$ such that $M+N$ is doubly stochastic.
\end{lem}

\begin{proof}
Let $M'$ be the $2n\times 2n$ doubly stochastic matrix such that
its upper left $n\times n$ sub-matrix is $M$. By Theorem~\ref{t double stoch},
\[ M'=\sum_{i=1}^S \theta_iP'_i\]
where $P'_i$ are $2n\times 2n$ permutation matrices and
$\sum_{i=1}^S \theta _i =1$. Suppose that $P_i$ is the upper left
$n\times n$ sub-permutation matrix of $P'_i$, then
\[ M=\sum_{i=1}^S \theta_iP_i\]
Let $Q_i$ be a $n\times n$ sub-permutation matrix such that
$P_i+Q_i$ is a permutation matrix, say $R_i$. Define
\[N=\sum_{i=1}^S \theta _iQ_i\]
thus \[M+N =\sum_{i=1}^S \theta _i R_i\] which is a doubly
stochastic matrix.
\end{proof}

\begin{lem}
\label{t sub stoch}
Let $M$ be an $n\times n$ matrix. If every row sum and column sum of
$M$ is 0 and $|M|$ is sub-doubly stochastic, then
\[M=\sum_{i=1}^S \theta _iP_i\] where $P_i$ are permutation matrices,
$\sum_{i=1}^S \theta _i=0$ and $\sum_{i=1}^S |\theta _i|=1$
\end{lem}

\begin{proof}
Let
\[A=\frac{|M|+M}{2}\]
\[B=\frac{|M|-M}{2}\]
so $A$ and $B$ are nonnegative, and $2A$ and $2B$ are sub-doubly
stochastic.
By Lemma~\ref{l M+N}, there exists a sub-doubly stochastic matrix $C$
such that $2(A+C)$ is a doubly stochastic. But $A$ and $B$ have
the same row sums and column sums, and hence $2(B+C)$ is also a doubly
stochastic. By applying Theorem~\ref{t double stoch}, we have
\[2(A+C)=\sum_{i=1}^m \lambda _iQ_i\]
\[2(B+C)=\sum_{i=1}^{m'} \lambda '_iQ'_i\]
where $Q_i$, $Q'_i$ are permutation matrices, and the $\lambda
_i$, $\lambda'_i$  are nonnegative numbers satisfying
$\sum_{i=1}^m \lambda_i =\sum_{i=1}^{m'} \lambda'_i =1$.
Then the result follows because
\[M=(A+C)-(B+C)=\sum_{i=1}^m \frac{\lambda _i}{2}Q_i -\sum_{i=1}^{m'} \frac{\lambda '_i}{2}Q'_i\]
\end{proof}

\begin{proof}[Proof of Proposition~\ref{t ek<=dk T disc}]
 From Lemma~\ref{t sub stoch}, we know that for each $k\geq1$ and
almost all $x_1,\dots,x_{k-1}$
\[T_k(x_1,\dots,x_{k-1})=\sum_{i_k=1}^{S_k} \theta _{k,i_k}(x_1,\dots,x_{k-1})\cdot P_{k,i_k}(x_1,\dots,x_{k-1})\]
where $P_{k,i_k}$ are permutation matrices, $\sum_{i=1}^{S_k}
\theta _{k,i_k}=0$, and $\sum_{i=1}^{S_k} |\theta _{k,i_k}|=1$.
Let
\begin{equation}
\label{h=Pd}
h_{k,i_k}(x_1,\dots,x_{k-1},\cdot)=[P_{k,i_k}(x_1,\dots,x_{k-1})]d_k(x_1,\dots,x_{k-1},\cdot).
\end{equation}
Then
\begin{eqnarray*}
e_k& =& \left[\sum_{i_k=1}^{S_k} \theta _{k,i_k}P_{k,i_k}\right]d_k\\
&=&\sum_{i_k=1}^{S_k} |\theta _{k,i_k}|\varepsilon _{k,i_k}h_{k,i_k}
\end{eqnarray*}
where $\varepsilon_{k,i_k}= \mbox{sgn}(\theta _{k,i_k})$.

Now we need to consider the probability space $\Omega_1 \times
\Omega_2$, where $\Omega_1 = \Omega_2 = [0,1]^{\mathbb{N}}$. We
consider all of the previous random variables as random variables
on this new probability space, depending only upon the first
coordinate $\omega_1$.  We define a filtration $(\mathcal G_k)$
where $\mathcal G_k = \mathcal F_k \otimes \mathcal L_{k+1}$.

We define a predictable sequence of random variables $(I_k)$ so that for
each $\omega_1 \in \Omega_1$, the random variable
$I_k(\omega_1,\cdot)$ takes the value $i$ with probability
$|\theta_{k,i}(\omega_1)|$.  Then we see that
\[ e_k = E[\varepsilon_{k,I_k} h_{k,I_k} |
   \mathcal L \otimes \{\emptyset,\Omega_2\}] .\]
Hence, since conditional expectation is a contraction on $L_p$
\[
\left\| \sum_{k=1}^n e_k \right\|_p \le \left\| \sum_{k=1}^n
\varepsilon_{k,I_k} h_{k,I_k} \right\|_p .
\]
Now we see that $(\varepsilon_{k,I_k})$ is a predictable sequence
bounded by $1$.  Hence by Burkholder's inequality, we see that
\[
\left\| \sum_{k=1}^n \varepsilon_{k,I_k} h_{k,I_k} \right\|_p \le
c_p \left\| \sum_{k=1}^n h_{k,I_k} \right\|_p .
\]
Next, observing (\ref{h=Pd}), since $P_{k,i_k}$ are permutation
matrices, for each $k\geq1$, $i_k=1,2,\dots,S_k$, $h_{k,i_k}$ is
just an $x_k$-rearrangement of $d_k$. that is
\[h_{k,i_k}(x_1,\dots,x_{k-1},j)=d_k(x_1,\dots,x_{k-1},\pi_{k,i_k}(j))\]
for some permutation $\pi _{k,i_k}$. Thus for any sequence
$(i_k)$ we have that
$(h_{k,i_k})$ and
$(d_k)$ are tangent sequences.
But then we see that $(h_{k,I_k})$ and $(d_k)$ are tangent
sequences.  Hence there
exists a positive constant $c_p$ such that
\[\left\|\sum_{k=1}^n h_{k,I_k}\right\|_p\leq
c_p \left\|\sum_{k=1}^n d_k\right\|_p . \]
The result follows.
\end{proof}

\begin{pro}
\label{t ek<=dk sharp disc}
Theorem~\ref{t ek<=dk sharp} holds in the case that
$(d_k)$ and $(e_k)$ are
adapted to the filtration $(\mathcal F_k)$ described above.
\end{pro}

This will follow immediately from the following well-known result
\cite{LT}.
\begin{thmofothers}
\label{t fsharp<=gsharp} $f=(f_1,f_2,\dots,f_N)$,
$g=(g_1,g_2,\dots,g_N)$ are $N$-dimensional real-valued vectors.
$f^\#=(f_1^\#,f_2^\#,\dots,f_N^\#)$ is the decreasing
rearrangement of $|f|=(|f_1|,|f_2|,\dots,|f_N|)$. Then
\[ \sum_{k=1}^n g_k^\#\leq \sum_{k=1}^n f_k^\# \]
for all $n=1,2,\dots,N$ if and only if there exists a matrix
$T=[a_{ij}]_{N \times N}$ such that $Tf=g$, $\sum_{i=1}^N
|a_{ij}|\leq 1$ and $\sum_{j=1}^N |a_{ij}|\leq 1$.
\end{thmofothers}

\section{The General Case}
The following theorem was proved by Crowe, Zweibel and Rosenbloom
\cite{CZR}.
\begin{thmofothers}
\label{t CZR}
Suppose $f$, $g$ are random variables on $[0,1]$, then for $1\leq
p\leq \infty $,
\[\left\|f^\# -g^\# \right\|_p\leq \|f -g\|_p\]
\end{thmofothers}

\begin{lem}\label{app} $1 \le p < \infty$. $(\Omega,\mathcal{F},P)$
is a probability space. Let $d(\omega,x)$,$e(\omega,x) \in
L_p(\Omega\times[0,1])$ be two random variables such that for
almost every $\omega \in \Omega$ that
$$
   \int_0^t (e(\omega,\cdot))^\#
   \le
   \int_0^t (d(\omega,\cdot))^\#
$$
and
$$ \int_0^1 d(\omega,\cdot) = \int_0^1 d(\omega,\cdot) = 0.
$$
Then given $\epsilon>0$, there exists a positive integer $N$ and
$d',e' \in L_p(\Omega\times[0,1])$ that are measurable with
respect to $\mathcal{F}\otimes\Sigma_N$ such that $\|d-d'\|_p,
\|e-e'\| \le \epsilon$,
$$
   \int_0^t (e'(\omega,\cdot))^\#
   \le
   \int_0^t (d'(\omega,\cdot))^\#
$$
and
$$ \int_0^1 d'(\omega,\cdot) = \int_0^1 d'(\omega,\cdot) = 0.
$$
\end{lem}

\begin{proof}
For every $\epsilon > 0$, pick $0 < \gamma < \min\{\epsilon /
[7(\|d\|_p \vee \|e\|_p)],1/3\}$ . Fix $\omega \in \Omega$, and
regard the functions as functions of only one variable $x$ on
$[0,1]$. Hence there exist simple functions
\[\bar{d}=\sum_{i=1}^S \bar{\alpha} _i\chi _{A_i}\]
\[\bar{e}=\sum_{i=1}^S \bar{\beta} _i\chi _{B_i}\]
such that
$$ \|\bar{d}- d\|_{L_p([0,1])}
   \le \gamma \| d \|_{L_1([0,1])}.$$
$$ \|\bar{e} - e\|_{L_p([0,1])}
   \le \gamma \| e \|_{L_1([0,1])}.$$
We may suppose without loss of generality that the sets $A_i$ and
$B_i$ are the sets of the form $[r_1,s_1)$, where the $r_i$ and
$s_i$ are rational numbers. Furthermore, we will suppose that
$A_{i_1} \cap A_{i_2} = \emptyset$ and $B_{i_1} \cap B_{i_2} =
\emptyset$ for $i_1 \ne i_2$.

Let $N_0=N_0(\omega)$ be the least common denominator of all these
rational numbers.  For each $\omega$, since $d(\omega,\cdot)^\#$,
$e(\omega,\cdot)^\#$ are Reimann integrable as a function of $x$,
there is a number $N_1=N_1(\omega)$ that is a multiple of $N_0$
and such that for all $n\ge N_1$ that
$$ \| E[d(\omega,\cdot)^\#|\Sigma_n] - d(\omega,\cdot)^\#\|_{L_p([0,1])}
   \le \gamma \| d \|_{L_1([0,1])}.$$
$$ \| E[e(\omega,\cdot)^\#|\Sigma_n] - e(\omega,\cdot)^\#\|_{L_p([0,1])}
   \le \gamma \| e \|_{L_1([0,1])}.$$
Now let $d_n = d \chi_{\{N_1(\omega) \le n\}}$ and $e_n = e
\chi_{\{N_1(\omega) \le n\}}$.  Then $d_n \to d$ and $e_n \to e$
in $L_p(\Omega\times[0,1])$.  So pick $N$ such that
$$\|d_N - d\|_{L_p(\Omega\times[0,1])} < \epsilon/7,$$
$$\|e_N - e\|_{L_p(\Omega\times[0,1])} < \epsilon/7.$$
For each fixed $\omega \in \{N_1(\omega) \le N\}$,
$[\frac{i-1}{N},\frac{i}{N})$ is either contained in some $A_j$ or
disjoint to all $A_j$. Let $\alpha_i=\bar{\alpha}_j$ if
$[\frac{i-1}{N},\frac{i}{N})\subset A_j$ for some $j$, and
$\alpha_i=0$ otherwise. Let $\chi
_i=\chi_{[\frac{i-1}{N},\frac{i}{N})}$. Thus
\begin{equation}
\label{approx-1} \left\|\sum_{i=1}^N \alpha _i\chi _i-d_N
\right\|_{L_p([0,1])} \leq \gamma \|d\|_{L_1([0,1])}
\end{equation}
and
\[
\left( \sum_{i=1}^N \alpha _i\chi _i\right)^\# = \sum_{i=1}^N
\varepsilon_{\sigma (i)}\alpha _{\sigma (i)}\chi _i
\]
for some permutation $\sigma $,  where $\varepsilon
_j=\mbox{sgn}(\alpha_j)$. By Theorem~\ref{t CZR},
\begin{equation}
\label{approx-3} \left\| \sum_{i=1}^N \varepsilon_{\sigma
(i)}\alpha _{\sigma (i)}\chi _i-d^\#_N \right\|_{L_p([0,1])} \le
\gamma \|d\|_{L_1([0,1])}
\end{equation}
and also the analogous statement holds for $e$.

Now if we set $$E[d^\#_N|\Sigma_N]=\sum_{i=1}^N \alpha''_i
\chi_i$$ then
\begin{equation}
\label{approx-4} \left\|\sum_{i=1}^N \alpha''_i\chi_i -d^\#_N
\right\|_{L_p([0,1])} \le \gamma \|d\|_{L_1([0,1])}
\end{equation}
Note that in this case that
\[
\int_0^t \sum_{i=1}^N \alpha''_i\chi _i=\int_0^t d^\#_N
\]
if $t=\frac{j}{N}$, $j=0,1,2,\dots,N$. Then by (\ref{approx-3})
and (\ref{approx-4}),
\[
\left\|\sum_{i=1}^N \alpha''_i\chi_i- \sum_{i=1}^N
\varepsilon_{\sigma (i)}\alpha _{\sigma (i)}\chi _i
\right\|_{L_p([0,1])} \le 2 \gamma \|d\|_{L_1([0,1])}
\]
By doing the reverse process of taking decreasing rearrangement of
$|\sum_{i=1}^N \alpha _i\chi _i|$, and setting
\[ \hat \alpha_i =
   \varepsilon_{\sigma^{-1}(i)}\alpha''_{\sigma^{-1}(i)} \]
we have
\begin{equation}
\label{approx-5} \left\|\sum_{i=1}^N \hat \alpha_i \chi_i -
       \sum_{i=1}^N \alpha _i\chi _i \right\|_{L_p([0,1])} \le
       2 \gamma \| d\|_{L_1([0,1])}
\end{equation}
 From (\ref{approx-1}) and (\ref{approx-5}),
\[
\left\|\sum_{i=1}^N \hat \alpha_i \chi_i-d_N \right\|_{L_p([0,1])}
       \le 3 \gamma \|d\|_{L_1([0,1])}
\]
For $t=\frac{j}{N}$, $j=0,1,2,\dots,N$, it is clear that
\[\int_0^t \left( \sum_{i=1}^N \hat{\alpha }_i\chi _i \right)^\# =
\int_0^t \sum_{i=1}^N \alpha'' _i\chi _i=\int_0^t d^\#_N.\]
Furthermore, if we set
\[ \zeta = E\left[\sum_{i=1}^N \hat{\alpha }_i\chi_i\right] \]
then $$|\zeta| \le 3 \gamma\|d\|_{L_1([0,1])}.$$ We can also
perform this same construction for $e$, the analogues of $\hat
\alpha_i$ and $\zeta$ being $\hat\beta_i$ and $\eta$. Thus we see
that for $t=\frac{j}{N}$, $j=0,1,2,\dots,N$ that
\begin{eqnarray}
\label{int ineq seq}
& &\int_0^t \left(\sum_{i=1}^N \left(\hat{\alpha}_i-\zeta \right)\chi _i \right)^\# \\
& \leq &\int_0^t\left(\sum_{i=1}^N\left(|\hat{\alpha}_i|+|\zeta |\right)\chi _i \right)^\#\nonumber\\
& \leq &\int_0^t \left( \sum_{i=1}^N \hat{\alpha}_i\chi _i
\right)^\# +
3\gamma \|d\|_{L_1([0,1])} \cdot t\nonumber\\
& =& \int_0^t d^\#_N + 3\gamma \|d\|_{L_1([0,1])} \cdot t\nonumber\\
& \leq & (1+3\gamma)\int_0^t d^\#_N \nonumber
\end{eqnarray}
and similarly
\begin{eqnarray}
\label{int approx} \int_0^t
\left(\sum_{i=1}^N\left(\hat{\beta}_i-\eta \right)\chi _i
\right)^\#
&\geq & \int_0^t e^\#_N - 3\gamma\|e\|_{L_1([0,1])} \cdot t\\
&\geq & (1-3\gamma)\int_0^t e^\#_N \nonumber
\end{eqnarray}
Thus, we are ready to define $d'$ and $e'$. Let
\[d'= (1+3\gamma) \sum_{i=1}^N (\hat{\alpha }_i-\zeta )\chi_i \]
\[e'= (1-3\gamma) \sum_{i=1}^N (\hat{\beta }_i-\eta )\chi _i\]
It is clear that $E[d']=E[e']=0 $. Combining (\ref{int ineq seq})
and (\ref{int approx}), we have for $t=\frac{j}{N}$,
$j=0,1,2,\dots,N$
\[ \int_0^t (e')^\# \leq \int_0^t e^\# _N =\int_0^t d^\# _N\leq \int_0^t (d')^\# . \]
But then by linear interpolation, this follows for all $t \in
[0,1]$. Now an easy argument shows that
\[ \|d'-d\|_{L_p(\Omega \times [0,1])}  \leq 6 \gamma\|d\|_{L_1(\Omega \times
[0,1])}+\epsilon/7\]
\[ \|e'-e\|_{L_p(\Omega \times [0,1])}  \leq 6 \gamma\|e\|_{L_1(\Omega \times
[0,1])}+\epsilon/7\] and we are done.
\end{proof}

\bigskip

\begin{proof}[Proof of Theorem \ref{t ek<=dk sharp}]
For each $1\leq k \leq n $, apply Lemma \ref{app}, there exists an
integer $N_k$ and functions $d'_k$, $e'_k$ satisfying
$\|d'_k-d_k\|_p$, $\|e'_k-e_k\|_p \leq \epsilon$ such that
$(d'_k)$ and $(e'_k)$ are adapted to $(\mathcal{L}_{k-1}\otimes
\Sigma_N)$, where $N$ is the least common multiple of $N_k$, keep
the martingale property, and
\[\int_0^t (e'_k(x_1,\dots,x_{k-1},\cdot))^{\#}(s)ds \leq \int_0^t  (d'_k(x_1,\dots,x_{k-1},\cdot))^{\#}(s)ds\]
for all $t\in [0,1]$. By Proposition~\ref{t ek<=dk sharp disc},
there exist a positive constant $c_p$ such that
\[ \left\|\sum_{k=1}^n e'_k \right\|_p
   \leq c_p \left\| \sum_{k=1}^n d'_k \right\|_p\]
 $\| d_k - d'_k \|_p \to 0$ and
$\| e_k - e'_k \|_p \to 0$ as $\epsilon \to 0$.
The result follows.
\end{proof}

\begin{proof}[Proof of Theorem~\ref{t ek<=dk T}]
If $f$ is a random variable on $(\Omega ,\mathcal{F},P)$, $1\leq
p<\infty$, $0\leq t\leq 1$, we define the $K$-functional by
\[K(t,f;L_p,L_\infty)=\inf_{f_0+f_1=f}\{\|f_0\|_p+t\|f_1\|_\infty \}.\]
J. Peetre \cite{P} has shown that
\[K(t,f;L_1,L_\infty) = \int_0^t f^{\#} (s)ds .\]
Hence it follows that if $T$ is an operator on both $L_1([0,1])$
and $L_\infty([0,1])$ with norm bounded by $1$, then for $t \ge 0$
\[ \int_0^t (Tf)^{\#}(s) ds \le \int_0^t f^{\#}(s) ds . \]
Thus the result follows from
Theorem~\ref{t ek<=dk sharp}.
\end{proof}

\begin{lem}
\label{l Mg<=Mf}
Let $f$ and $g$ be real-valued random variables on
$(\Omega ,\mathcal{F},P)$.  Then
\begin{equation}
\label{Mg<=Mf} E \left[ \lambda \vee |g| \right]\leq
E\left[\lambda \vee |f| \right]
\end{equation}
for all nonnegative number $\lambda$  if and only if
\[\int_0^t g^{\#}(s) ds\leq \int_0^t f^{\#}(s) ds \]
for all $t\in [0,1]$.
\end{lem}

\begin{proof}
Equation~(\ref{Mg<=Mf}) is equivalent to $E \left[ \lambda \vee
g^\#  \right] \leq E \left[\lambda \vee f^\# \right]$. For the
``if'' part, let
\[\alpha =\sup \left\{t:f^\# (t) \geq \lambda \right\}\]
\[\beta =\sup \left\{t:g^\# (t)\geq \lambda \right\}.\]
Then
\begin{eqnarray*}
E \left[ \lambda \vee f^\#\right] &=&
\int_0^\alpha f^\#  + (1-\alpha)\lambda \\
&=&
\int_0^\beta f^\# + (1-\beta)\lambda + \int_\beta^\alpha (f^\#-\lambda) \\
&\ge&
\int_0^\beta g^\# + (1-\beta)\lambda + \int_\beta^\alpha (f^\#-\lambda) \\
&=& E \left[ \lambda \vee g^\# \right] + \int_\beta^\alpha
(f^\#-\lambda) .
\end{eqnarray*}
If $\alpha \le \beta$, then for all $x \in (\alpha,\beta)$ we have
$f^\#(x) \le \lambda$, and if $\beta \le \alpha$, then for all $x
\in (\beta,\alpha)$ we have $f^\#(x) \ge \lambda$.  Either way, we
see that $\int_\beta^\alpha (f^\#-\lambda) \ge 0$, and the result
follows.

To show the ``only if'', for any $\alpha \in [0,1]$, let
\[\lambda=f^\#(\alpha)\]
\[\beta =\inf \left\{t:g^\# (t)\geq \lambda \right\}.\]
Then
\begin{eqnarray*}
\int_0^\alpha g^\#
&=&
\int_0^\beta g^\# + \int_\beta^\alpha (g^\#-\lambda) + \lambda(1-\beta) + \lambda(\alpha-1) \\
&=&
E \left[ \lambda \vee g^\# \right] + \lambda(\alpha-1) + \int_\beta^\alpha (g^\#-\lambda) \\
&\le&
E \left[ \lambda \vee f^\# \right] + \lambda(\alpha-1) + \int_\beta^\alpha (g^\#-\lambda) \\
&=& \int_0^\alpha f^\# + \int_\beta^\alpha (g^\#-\lambda) .
\end{eqnarray*}
Arguing as above, we see that $\int_\beta^\alpha (g^\#-\lambda)
\le 0$, and again the result follows.
\end{proof}

Given a random variable $f$ and a sigma field $\mathcal{G}$, we
will say that $f$ is nowhere constant with respect to
$\mathcal{G}$ if $P(f=g)=0$ for every $\mathcal{G}$ measurable
function $g$. The following theorem \cite{M2} shows a concrete
representation of a sequence of random variables.

\begin{thmofothers}
\label{t concrete} Let $(f_n)$ be a sequence of random variables
takeing values in a separable sigma filed $(S,\mathcal{S})$. Then
there exists a sequence of measurable functions
$(g_n:[0,1]^{n}\rightarrow S)$ that has the same law as $(f_n)$.
If further we have that $f_{n+1}$ is nowhere constant with respect
to $\sigma(f_1,\dots,f_n)$ for all $n \geq 0$, then we may suppose
that $\sigma(g_1,\dots,g_n)=\mathcal{L}_n$ for all $n \geq 0$.
\end{thmofothers}

\begin{proof}[Proof of Theorem~\ref{t ek<=dk sk}]
We will prove this theorem under the assumption
(\ref{skekFk-1<=skdkFk-1}). Consider the map
$D_k=(d_k,e_k,f_k):\Omega \times [0,1]^\mathbb{N} \rightarrow
\mathbb{R}^3$ by $(\omega,(x_k))\mapsto
(d_k(\omega),e_k(\omega),x_k)$. It is clear that $D_k$ is nowhere
constant with respect to $\sigma(D_1,\dots,D_{k-1})$. Apply the
previous theorem to get
$\widetilde{D}_k=(\widetilde{d}_k,\widetilde{e}_k,\widetilde{f}_k):[0,1]^k
\rightarrow \mathbb{R}^3 $ such that $(\widetilde{D}_k)$ has the
same law as $(D_k)$ and
$\sigma(\widetilde{D}_1,\dots,\widetilde{D}_k)=\mathcal{L}_k$.

Next, we show that for almost every $x_1, \dots,x_{k-1}$ and
$\lambda \ge 0$ that
\[ \int_0^1 \lambda \vee |\widetilde e_k(x_1,\dots,x_k)| \, dx_k
   \le
   \int_0^1 \lambda \vee |\widetilde d_k(x_1,\dots,x_k)| \, dx_k \]
which will follow from showing that for any
bounded non-negative measurable function
$\phi_k : [0,1]^{k-1} \rightarrow [0,\infty)$ that
\[ E[\phi_k \vee |\widetilde{e}_k|]
\leq E[\phi_k \vee |\widetilde{d}_k|] .\] But then there exists a
bounded Borel measurable function $\theta_k : \mathbb{R}^{3(k-1)}
\rightarrow [0,\infty)$ such that
$\phi=\theta(\widetilde{D}_1,\dots,\widetilde{D}_{k-1})$ almost
everywhere in $[0,1]^{k-1}$. Thus
\begin{eqnarray*}
\int_{[0,1]^k} \phi_k \vee |\widetilde{e}_k| &=&\int_{[0,1]^k}
\theta(\widetilde{D}_1,\dots,\widetilde{D}_{k-1})\vee
   |\widetilde{e}_k| \\
&=& E[ \theta(D_1,\dots,D_{k-1}) \vee |e_k|] \\
&\le& E[ \theta(D_1,\dots,D_{k-1}) \vee |d_k|] \\
&=&\int_{[0,1]^k}
\theta(\widetilde{D}_1,\dots,\widetilde{D}_{k-1})\vee
   |\widetilde{d}_k| \\
&=&
\int_{[0,1]^k} \phi_k \vee |\widetilde{d}_k|
\end{eqnarray*}
Also to show that
$E[\widetilde{d}_k|\mathcal{L}_{k-1}]=E[\widetilde{e}_k|\mathcal{L}_{k-1}]=0$,
it is sufficient to show that for any bounded measurable function
$\phi_k : [0,1]^{k-1} \rightarrow \mathbb{R}$ that $E[\phi_k
\widetilde{d}_k]=E[\phi_k \widetilde{e}_k]=0$.  Thus follows by a
very similar argument to that above.

The result then
follows from Lemma~\ref{l Mg<=Mf} and Theorem~\ref{t ek<=dk sharp}.
\end{proof}

{\bf Acknowledgments.} We would like to mention the help of Jim
Reeds and David Boyd in obtaining the argument about doubly
stochastic matrices and other useful remarks.  Also, we would like to
express our thanks
to Mitch Taibleson for bringing this problem to our attention.

\bibliographystyle{amsplain}

\end{document}